\documentstyle[12pt]{article}
\begin{document} 
\input{mssymb}
\setlength{\textwidth}{6in}
\setlength{\textheight}{8in}
\setlength{\topmargin}{0.3in}
\setlength{\headheight}{0in}
\setlength{\headsep}{0in}
\setlength{\parskip}{0pt}

\newcounter{probnumt}
\partopsep 0pt
\settowidth{\leftmargin}{1.}\addtolength{\leftmargin}{\labelsep}
\newenvironment{probt}{\begin{quote}\begin{enumerate} 
\setcounter{enumi}
{\value{probnumt}}}%
{  \setcounter{probnumt}
{\value{enumi}}\end{enumerate}\end{quote}}

\newcommand{\RR}{{\Bbb R}}
\newcommand{\CC}{{\Bbb C}}
\newcommand{\KK}{{\Bbb K}}
\renewcommand{\limsup}{\overline {\lim}\,}
\renewcommand{\liminf}{\underline {\lim}\,}
\newcommand{\two}{\mbox {2-polynomial}}
\newtheorem{thm}{Theorem}
\newtheorem{lem}[thm]{Lemma}
\newtheorem{cor}[thm]{Corollary}
\newtheorem{pro}{Question}

\newenvironment{proof}{\medskip\par\noindent{\bf Proof.}}{\hfill
$\Box$ \medskip \par}

\setcounter{page}{1}

\title{On the extension of 2-polynomials}

\author { Pei-Kee Lin\\
Department of Mathematics \\
University of Memphis\\
Memphis TN, 38152\\
USA}

\maketitle
%\footnotetext {1991 Mathematics Subject Classification: 46A22}
%\footnotetext{Key words and phrases: $\two$, bilinear functional,
%inner product space}

%\begin{abstract}
%{Let $X$ be a three dimensional real Banach space.
%Ben\'itez and Otero \cite {BeO} showed that
%if the unit ball of $X$ is 
%is an  intersection of two ellipsoids, then
%every $\two$ defined in a linear subspace of $X$ can
%be extended to $X$ preserving the norm.  In this article,
%we extend this result to any finite dimensional Banach space.}
%\end{abstract}

Let $X$ be a normed linear space over $\KK$ ($\RR$ or $\CC$). A function
$P:X \to \KK$ is said to be a 2-{\em polynomial} if there is a bilinear
functional $\Pi:X \times X \to \KK$ such that $P(x)=\Pi(x,x)$ for
every $x \in X$. 
The norm of $P$ is defined by
\[\|P\|=\sup \{|P(x)|:\|x\|=1\}.\]
It is known that if $X$ is an inner product
space, then every $\two$ defined in a linear subspace of $X$ can
be extended to $X$ preserving the norm.  On the other hand, there
is a $\two$ $P$ defined in a two dimensional subspace of
$\ell_\infty^3$ such that every extension of $P$ to
$\ell_\infty^3$ has norm greater than $\|P\|$, (see \cite {AB,
BeO}).  Recently, Ben\'itez and Otero \cite {BeO} showed that
if $X$ is a three dimensional real Banach space $X$ such that the
unit ball of $X$ is an  intersection of two ellipsoids, then
every $\two$ defined in a linear subspace of $X$ can
be extended to $X$ preserving the norm.  It is natural to ask
\cite {BeO}

\begin{pro} Suppose $X$ is a norm space such that the
unit ball of $X$
is an intersection of two ellipsoids.  Can
every $\two$ defined in a linear subspace of $X$
be extended to $X$ preserving the norm?
\end{pro}

In this article, we show the answer is affirmative when $X$ is 
a finite dimensional space.  

First, we
recall the following  result in \cite {BeO}.

\begin{lem}\label{A} 
If $P$ and $Q$ are $\two$s in $X=\RR^2$ such that
\[ 0 \leq \sup\{P(x), Q(x)\}, \hspace{.5in} (x \in X)\]
then there exists $0 \leq \alpha \leq 1$ such that
\[0 \leq \alpha P(x)+ (1-\alpha) Q(x), \hspace {.5in} (x \in X).
\]
\end{lem}

%\begin{lem}\label{B} Suppose $X$ be a two dimensional real space.
%Let $P$ be an indefinite $\two$ on $X$, and $\Pi$ be the
%symmetric bilinear functional associated to $P$, and $x_1,x_2 \in
%S=\{\|x\|=1: x \in X\}$ be such that
%\[P(x_1)=\sup\{P(x):x \in S\}, \hspace{.5in}
%P(x_2)=\inf\{P(x):x \in S\}.\]
%Then for any $x \in S$
%\[ P(x_i)P(x) \leq \Pi^2(x_i,x) \leq P^2(x_i), \hspace{.5in}
%(i=1,2).
%\]
%\end{lem}
 
\begin{lem}\label{C}  Suppose $X$ is a real norm space, and suppose
$P_1, P_2$ are two positive semidefinite 
$\two$s on $X$. If $Q$ is a $\two$ such that 
\[|Q(x)| \leq \max \{P_1(x), P_2(x)\} \hspace {.5in} (x \in X),\]
then there are $0 \leq \alpha, \beta \leq 1$ such that
for every $x \in X$,
\[-((\beta P_1(x) +(1-\beta)P_2(x))\leq Q(x)
\leq \alpha P_1(x) +(1-\alpha)P_2(x).\]
\end{lem}

\begin{proof}  For every $x \in S$ with $Q(x) \geq 0$ ($Q(x) \leq 0$),
let 
\[A_x=\{ \gamma \in [0,1]: \gamma P_1(x)+(1-\gamma) P_2(x) \geq Q(x) \}
\]
\[\bigl (B_x =\{\gamma \in [0,1]:-(\gamma P_1(x)+(1-\gamma)
P_2(x)) \leq Q(x) \}.\bigr)\] 
To prove this lemma, it is enough to show that
\[\cap_{\{x \in S :
 Q(x) \geq 0\}} A_x\ne \emptyset, \mbox { and }
\cap_{\{x \in S 
: Q(x) \leq 0\}} B_x \ne \emptyset.\] 
For any $x \in S$ with $Q(x)\geq 0$ (respectively,  $Q(x) \leq
0$), $A_x$ (respectively, $B_x$) is a nonempty closed 
subinterval of $[0,1]$. So we only
need to show that for any $x,y \in S$, if $\min \{Q(x),Q(y)\} \geq
0$, (respectively, $\max\{Q(x),Q(y)\} \leq 0$), then
$A_x \cap A_y \ne
\emptyset$, (respectively, $B_x \cap B_y \ne \emptyset$).  

Suppose that $Q(x) \geq 0$ and $Q(y) \geq 0$.  
For any $z \in \mbox{span} \, \{x,y\}$, $R_1$ and $R_2$
are defined by
\begin{eqnarray*}
R_1(z)=&P_1(z)-Q(z) \\
R_2(z)=&P_2(z)-Q(z).
\end{eqnarray*}
By Lemma~\ref{A}, there is a $\gamma$, $0 \leq \gamma \leq 1$
such that
\[0 \leq \gamma R_1(z) +(1-\gamma)R_2. \]
This implies $\gamma \in A_x \cap A_y$.
Similarly, if $Q(x) \leq 0$ and $Q(y) \leq 0$, then there exists
$0 \leq \gamma \leq 1$ such that $\gamma \in B_x \cap B_y$.
We proved our lemma. \end{proof}

\begin{thm}  Let $\Pi_1$ and $\Pi_2$ be two inner
products on $\RR ^n$ and let $X=\RR^n$ be the space with the
norm
\[\|x\|=\sqrt {\max\{\Pi_1(x,x), \Pi_2(x,x)\}}.\]
Then every $\two$ defined in a subspace of $X$ can be extended to
$X$ preserving the norm.
\end{thm}
\begin{proof}  Let $P$ be any $\two$ on a subspace $Y$ of $X$
and let
\begin{eqnarray*}
\|x\|_1&=\sqrt{\Pi_1(x,x)} \\
\|x\|_2&=\sqrt{\Pi_2(x,x)}.\end{eqnarray*}
Without loss of generality, we may assume $\|P\|=1$
and $Y$ is a hyperplane of $X$.
So for any $x \in Y$,
\[
|P(x)| \leq \|x\|^2 \leq \max\{\|x\|_1^2, \|x\|_2^2\}.\]
By Lemma~\ref{C}, there exist $0 \leq \alpha, \beta \leq 1$ such
that for every $x \in Y$,
\[-(\beta \|x\|_1^2 +(1-\beta)\|x\|_2^2)\leq P(x)
\leq \alpha \|x\|_1^2 +(1-\alpha)\|x\|_2^2.\]
Replacing $\|x\|_1$ (respectively, $\|x\|_2$) by $(\alpha
\Pi_1(x,x)+(1-\alpha)\Pi_2(x,x))^{1/2}$ (respectively, $(\beta
\Pi_1(x,x)+(1-\beta)\Pi_2(x,x))^{1/2}$),
we may assume that if $x \in Y$, then
$$-\|x\|_2^2 \leq P(x) \leq \|x\|_1^2.$$

Let $\Pi$ be the symmetric bilinear functional associated to $P$.
Then there is a (bounded) symmetric
operator $T_1$ (respectively, $T_2$) on $(Y,\|\cdot\|_1)$
(respectively, $(Y,\|\cdot\|_2)$) such that for any $y_1,y_2 \in
Y$ 
\[ \Pi(y_1,y_2)=\Pi_1(y_1,T_1(y_2))=\Pi_2(y_1,T_2(y_2)).
\]
Since $T_1$ and $T_2$ are symmetric, they are diagonalizable. 
It is known that
\begin{probt}
\item every eigenvalue of $T_1$ (respectively, $T_2$) is real;
\item if $x_1$ and $x_2$ are two eigenvectors of $T_1$
(respectively, $T_2$) associated with two distinct eigenvalues,
then
\[\Pi (x_1,x_2)=0.\]
\end{probt}
Let $Y_1$ (respectively, $Y_3$) be the subspace  spanned by
all eigenvectors of $T_1$ (respectively, $T_2$) associated with non-negative
eigenvalues, and  $Y_2$ (respectively, $Y_4$) be the subspace  spanned by
all eigenvectors of $T_1$ (respectively, $T_2$) associated with negative
eigenvalues.  Then
$Y_1,Y_2, Y_3,Y_4$ satisfy the following conditions:
\begin{probt}
\item \label {a} $Y=Y_1 \oplus Y_2 =Y_3\oplus Y_4$;
\item \label {b}   $T_1(Y_1) \subseteq Y_1$, $T_1(Y_2) \subseteq
Y_2$, $T_2(Y_3) \subseteq Y_3$, and $T_2(Y_4) \subseteq Y_4$;
\item \label {c}  for any $y_i \in Y_i\setminus\{0\}$,
\[\begin{array}
{rl}
\Pi_1(y_1,T_1 (y_1)) &\geq 0 >\Pi_1(y_2,T_1 (y_1)),\\
\Pi_2(y_3,T_1 (y_3)) &\geq 0 >\Pi_2(y_4,T_2 (y_4)),\\
\Pi_1(y_1,y_2)&=0=\Pi_2(y_3,y_4).
\end{array}\]\end{probt}
We claim that $Y_1 \cap Y_4=\{0\}$.
Suppose it is not true. Let  $y \in (Y_1 \cap Y_4)\setminus
\{0\}$. Then 
\[ 0 > \Pi_2(y,T_2 (y))=\Pi(y,y)=\Pi_1(y,T_1(y)) \geq 0.\]
We get a contradiction.  Similarly, $Y_2 \cap Y_3=\{0\}$.
Hence, we have
\begin{probt}
\item \label {d} $\mbox {dim}\, (Y_1)=\mbox {dim}\,
(Y_3)$,   $\mbox {dim}\, (Y_2)=\mbox {dim}\,
(Y_4)$, 
and 
\[Y_1 \oplus Y_4=Y=Y_2\oplus Y_3.\]
\end{probt}
Let
\begin{eqnarray*}
M_1&=\{z\in X: \Pi_1(z,x)=0 \mbox { for all $x \in Y_1$}\},\\
M_2&=\{z\in X: \Pi_2(z,x)=0 \mbox { for all $x \in Y_4$}\}.
\end{eqnarray*}
By (\ref{c}) and (\ref {d}), 
$\mbox{dim}\, (M_1)=\mbox {dim} \, (Y_2)+1$ and 
$\mbox{dim}\, (M_2)=\mbox {dim} \, (Y_3)+1$.
This implies
there is a non-zero vector $z \in M_1 \cap
M_2$.  Let $\phi$ be any non-zero linear functional on $X$ such that
$\mbox {ker}\, \phi=Y$.           For any $x \in X$, define
\[\tilde P(x)=P\left (x-\frac {\phi(x)}{\phi(z)} z\right ).\]
We claim that if 
$0 <\tilde P(x)$, then $\tilde P(x) \leq \|x\|^2$.

Case 1.  $x -\frac{\phi(x)}{\phi(z)} z \in Y_1$.  Since $z \in
M_1$,          we have $\Pi_1(x-\frac {\phi(x)}{\phi(z)}z,z)=
0$.  So
\[ \tilde P(x)=P\left(x -\frac{\phi(x)}{\phi(z)} z\right) \leq
\left \|x -
\frac{\phi(x)}{\phi(z)} z\right\|^2_1 \leq \|x\|^2_1.\]

Case 2.  $x -
\frac{\phi(x)}{\phi(z)} z \notin Y_1$.  Then there exist $y_1 \in
Y_1$ and $y_2\in Y_2$ such that $x -
\frac{\phi(x)}{\phi(z)} z =y_1+y_2$.  Note: $\Pi_1(y_1,y_2)=0
=\Pi_1(y_1,z)$. So
\[\begin{array}{rll}
\tilde P (x) & = P\left(x -
\frac{\phi(x)}{\phi(z)} z\right)
&=   \Pi(y_1+y_2,y_1+y_2)\\
&= P(y_1)+P(y_2) &\leq  P(y_1) \\
&\leq   \|y_1\|_1^2  &\leq \|x\|_1^2\\
&\leq\|x\|^2.
\end{array}\]
We  proved our claim.  Similarly, if $\tilde P(x) \leq 0$,
then $|\tilde P(x)| \leq \|x\|^2$.
Hence, $\tilde P$ is an extension of $P$ preserving the norm.
\end{proof}

In \cite {BeO} (Lemma 2 and Proposition 2), 
Ben\'itez and Otero proved that the problem of extension 
preserving the norm can be reduced to the real case.
Hence, we have the following theorem.

\begin{thm}  Let $\Pi_1$ and $\Pi_2$ be two inner
products on $\CC^n$ and let $X=\CC^n$ be the space with the
norm
\[\|x\|=\sqrt {\max\{\Pi_1(x,x), \Pi_2(x,x)\}}.\]
Then every $\two$ defined in a subspace of $X$ can be extended to
$X$ preserving the norm.
\end{thm}

\begin{center} Acknowledgements \end{center}

The author are thankful to the referee for his
comments.

\end{document}